\documentclass{autart}
\usepackage{amsmath,amssymb}
\usepackage{array}
\usepackage{mdwmath}
\usepackage{mdwtab}
\usepackage[mathscr]{eucal}
\usepackage{graphicx}
\usepackage{times}
\usepackage{pifont}
\usepackage{srcltx}
\usepackage{subfigure}
\usepackage{color}
\usepackage{epsfig}
\usepackage{enumerate}

\newcommand{\ov}[1]{\overline{#1}}

\newcommand{\br}{\mathbb{R}}

\newcommand{\sgn}{\text{\rm sgn}}
\newcommand{\url}[1]{\texttt{#1}}
\renewcommand{\Box}{\bullet}

\newcommand{\bul}{\bullet}

\newcommand{\sclp}[2]{\left\langle #1;  #2 \right\rangle}

\newcommand{\blr}{\boldsymbol{r}}
\renewcommand{\thefootnote}{\fnsymbol{footnote}}
\newtheorem{remark}{Remark}

\newtheorem{lemma}{Lemma}
\newtheorem{theorem}{Theorem}
\newtheorem{assumption}{Assumption}
\newcommand{\epf}{$\quad \bullet$}

\begin{document}
\begin{frontmatter}
\title{Proofs of the Technical Results Justifying an Algorithm for Collision Avoidance in Dynamic Environments with Moving and Deforming Obstacles}
\author[auth1]{Chao Wang}\ead{mch.hoy@gmail.com},
\author[auth2]{Alexey S. Matveev}\ead{almat1712@yahoo.com},
\author[auth1]{Andrey V. Savkin}\ead{a.savkin@unsw.edu.au}

\address[auth1]{School of Electrical Engineering and
Telecommunications, The University of New South Wales, Sydney 2052,
Australia}
\address[auth2]{Department of Mathematics and Mechanics, Saint Petersburg
University, Universitetskii 28, Petrodvoretz, St.Petersburg, 198504,
Russia}
\begin {abstract}
This text presents the proofs of the technical facts underlying theoretical justification of the convergence and performance of the novel algorithm
for reactive navigation of differential drive wheeled robots in dynamic uncertain environments. The algorithm restricts neither the natures nor the motions of the obstacles, they need not be rigid but conversely may deform. It does not consume data about the velocities, shapes, sizes, or orientations of the obstacles, and does not need a map of the environment or recognition of individual obstacles.
The only information about the scene is the current distance to the nearest obstacle.
\end{abstract}
\end{frontmatter}
\section{Introduction}
A key component of safe navigation is avoidance of collisions with en-route obstacles.
This problem has attracted an enormous attention in robotics research. The relevant algorithms can be generally classified into global and local path planners \cite{LZL07}.
\par
Global planners typically build a more or less comprehensive model of the environment to find the best complete trajectory \cite{Latom91}. By and large, they are computationally expensive; NP-hardness was established for even the simplest problems of dynamic motion planning \cite{Canny88}. This seriously troubles their real-time implementation.
Data incompleteness and erroneousness, typical for onboard perception, may cause a noticeable deterioration in the overall performance of global planners.
\par
Conversely, local planners iteratively re-compute a short-horizon path based on sensory data about a close vicinity of the robot. A short computation time typical for these planners
creates a potential for their use in real-time guidance systems.
However most of the related techniques either in fact treat the obstacles as static or assume a deterministic knowledge about the obstacle velocity and a moderate rate of its change; the first group is exemplified by the dynamic window \cite{SeMaPe05,FBTh97}, the curvature velocity \cite{Simm96}, and the lane curvature \cite{NaSi98} approaches, whereas velocity obstacles
\cite{FS98}, collision cones \cite{ChGh98}, or inevitable collision states \cite{FraAs03,OwMo06} methods provide examples of the second kind.
\par
In the marginal case where the planning horizon concentrates into a point, local planner acts as a reactive controller: it provides a reflex-like control response to the current observation. Examples include artificial potential approach \cite{LinHuLa06,FeRu07} and kinematic control based on polar coordinates and Lyapunov-like analysis \cite{ChQuPoFa10}. In this area, the obstacles were mostly viewed as rigid bodies of the simplest shapes \cite{FeRu07,ChQuPoFa10,LinHuLa06,MaKa07}, the sensory data were assumed to be enough to provide access to the locations of the obstacle characteristic points concerned with its global geometry (e.g., the disc center \cite{FeRu07,ChQuPoFa10} or angularly most distant polygon vertex \cite{MaKa07}) and to the full velocity of the obstacle \cite{FeRu07,MaKa07,ChQuPoFa10}. Furthermore, rigorous justification of the global convergence of the proposed algorithm was rarely encountered.
\par
In this paper, the problem of reactive navigation is addressed for a standard differential drive wheeled robot. The objective is to drive the robot to the assigned point-wise target through the obstacle-free part of the scene. The robot may have a limited sensor range, may not distinguish between various  points of the obstacle, and so may be unable to estimate many of its parameters, like size, center, edge, full velocity, etc. It has access only to the distance to the nearest obstacle point within the sensor range, the time derivative of this reading, and the angle-of-sight at the target.
Unlike the overwhelming majority of the previous works, the obstacles are not assumed to be rigid or even solid: they are continuums of arbitrary and time-varying shapes undergoing general tranformatons, including rotations and deformations.
\par
The extended introduction and discussion of the proposed control algorithm are given in a paper submitted by the authors to the IEEE Transactions on Control Systems Technology. This algorithm was extensively tested via both computer simulations and real-world experiments with both research wheeled robots and an intelligent wheelchair. The results of both simulation tests and experiments are also presented in that paper.
This text basically contains the proofs of the technical facts underlying theoretical justification of the convergence and performance of the proposed algorithm. These proofs were not included in that paper due to the length limitations. To make the current text logically consistent and self-contained, we supply the reader with relevant notations and complete statements of both the problem and theorems to be proved in this text.
\section{System Model and Problem Statement}\label{s2}
We consider a planar differential drive wheeled robot (DDWR) with two independently actuated driving wheels mounted on the same axle and castor wheels. DDWR is controlled by the angular velocities $\omega_l$ and $\omega_r$ of the left and right driving wheels, respectively.
These velocities are limited by a common constant $\Omega \geq |\omega_l|, |\omega_r|$; the
driving wheels roll without sliding.
\begin{figure}[ht]
\centering
\scalebox{0.3}{\includegraphics{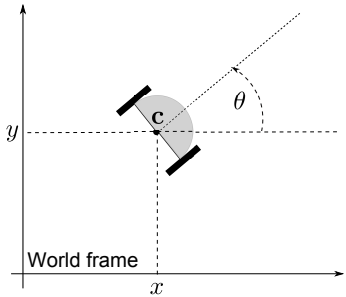}}
\caption{Differential drive wheeled robot}
\label{fig.ddwr}
\end{figure}
\par
The kinematics of such vehicles are classically described by the following equations:
\begin{equation}
\label{1}
\begin{array}{l}
\dot{x} = v \cos \theta,
\\
\dot{y} = v \sin \theta,
\\
\dot{\theta} = u,
\end{array}\;
\begin{array}{l}
v = \frac{v_l+v_r}{2},
\\
u = \frac{v_r-v_l}{2L},
\\
v_i = R_w \omega_i,
\end{array}
\;
\begin{array}{l}
x(0) = x_0
\\
y(0) = y_0
\\
\theta(0) = \theta_0
\end{array}.
\end{equation}
Here $x,y$ and $\theta$ are introduced in Fig.~\ref{fig.ddwr} and stand for the position and orientation of DDWR, $v$ is its longitudinal speed,
$R_w$ is the radius of the driving wheels, $2L$ is the length of the axle, and
$\omega_i = \omega_i(t) \in [-\Omega,\Omega], i=l,r$. To simplify the matters, we treat $v$ and $u$ as control variables. They uniquely determine the rotational velocities $\omega_r=(v+Lu)/R_w, \omega_l=(v-Lu)/R_w $ and obey the bound:
\begin{equation}
\label{control.constr}
|v|+L|u| \leq V:= R_w \Omega.
\end{equation}
It follows that for given $v \in (-V,V)$, the turning radius of DDWR is no less than
\begin{equation}
\label{mini}
R = \frac{L|v|}{V-|v|}.
\end{equation}
\par
The workspace of DDWR contains a moving obstacle $D(t) \subset \br^2$, which need not be a rigid body but conversely, may undergo arbitrary motions, including deformations. Its position and shape are not known in advance; the only available information is the current distance $d(t) := \text{\bf dist}_{D(t)}[r(t)]$ to the obstacle and the rate $\dot d(t)$ at which this reading evolves over time. Here $r:=[x, y]^{\top}$ is the vector of the coordinates of DDWR and
\begin{equation}
\label{dist} \text{\bf dist}_{D}(r) :=
\min_{r^\prime \in D}
\|r-r^\prime\|,
\end{equation}
where $\|\cdot\|$ denotes the standard Euclidean vector norm and the minimum is achieved if $D$ is closed.
\par
Finally, there is a steady point-wise target $\bf{T}$ in the plane and DDWR has constantl access to the heading $h(t)$ towards the target.
The objective is to guide DDWR through the obstacle-free part of the  plane and reach the target $\bf{T}$ at a certain time $t_f>0$:
	\begin{equation*}
				r(t_f) = \text{\bf T};
		\qquad
            r(t) \not\in D(t)\quad \forall t \in [0,t_f].
			\end{equation*}
Moreover, the distance to the obstacle should constantly exceed a given safety margin $d_{\text{safe}}>0$:
	\begin{equation}
	\label{dist}
		\text{\bf dist}_{D(t)}\left[ r(t) \right]\geq d_{\text{safe}}\quad \forall t \in [0,t_f].
	\end{equation}
\section{Autonomous Navigation Algorithm}\label{s3}
Now we present a summary of the proposed reactive navigation algorithm. The algorithm combines obstacle avoidance behavior, which is activated in a close proximity of the obstacle, with motions towards the target in a straight line when there is no threat of collision.
\par
We employ the following obstacle avoidance strategy:
\begin{gather}
\label{avoid}
\left.
 \begin{array}{l}
	u(t) = \frac{V - v(t)}{L} \cdot \sgn \big\{ \dot{d}(t) + \chi[d(t)-d_0]\big\}
\\
v(t) = \Upsilon[d(t)]
\end{array} \right|,
\;
 	\text{where}
 	\\
 	\label{chi}
	\chi(z):=
	\begin{cases} \gamma z & \text{if} \; |z|\leq \delta
	\\
	v_{\ast} \sgn (z)& \text{if} \; |z| > \delta
	\end{cases} \qquad (v_{\ast} := \gamma\delta )
\end{gather}
is the linear function with saturation, the smooth function $\Upsilon(\cdot) : [0,\infty) \to (0,V)$ determines the speed $v$ of DDWR on the basis of the current distance $D$ to the obstacle,
whereas $\gamma>0, \delta >0$, and $d_0 > d_{\text{safe}}$ are the controller parameters. Here  $d_0$ has the sense of the desired distance to the obstacle when bypassing it. The function $\Upsilon(\cdot)$ varies between two speeds $v_0$ and $v_{\text{cr}}$, i.e., $\Upsilon(d) = v_0 \; \forall d \leq d_0^\Upsilon, \Upsilon(d) = v_{\text{cr}} \; \forall d \geq d_{\text{cr}} > d_0^\Upsilon$. The speed $v_0 \in (0,V)$ is used when bypassing obstacles, so $d_0<d_0^\Upsilon$; the larger $v_{\text{cr}} >v_0$ \textquoteleft cruise\textquoteright\, speed is employed where there is no collision threat. We assume that $v_{\text{cr}} <V$ to leave the vehicle a certain level of maneuverability in the "cruise" regime.
\par
The proposed obstacle avoidance strategy belongs to the class of sliding mode control algorithms; see e.g. \cite{UTK92}.
The intuition behind this strategy is that in the sliding mode, the equation $\dot{d}+\chi(d-d_0)=0$ of the sliding surface is satisfied,
according to which the vehicle is steered towards the desired distance $d_0$ to the obstacle: $d \to d_0$. For this to take effect, the sliding mode maneuver should be at least realistic. Since the derivative $\dot{d}$ does not exceed $|\dot{d}| \leq v_r$ the relative speed $v_r$ of the vehicle with respect to the obstacle, this means that in (\ref{chi}), the saturation level $v_\ast$ should not exceed this speed. This can be achieved by proper tuning of the controller parameters $\gamma$ and $\delta$ based on available estimates of the speed of the obstacle.
If initially the vehicle is not on the sliding surface, the control law \eqref{avoid} quickly drives it to this surface after a short initial turn, see Section~\ref{sec.proof} for details. So sliding motion is the main part of the obstacle avoidance maneuver.
\par
\begin{figure}[h]
\centering
\subfigure[]{\scalebox{0.28}{\includegraphics{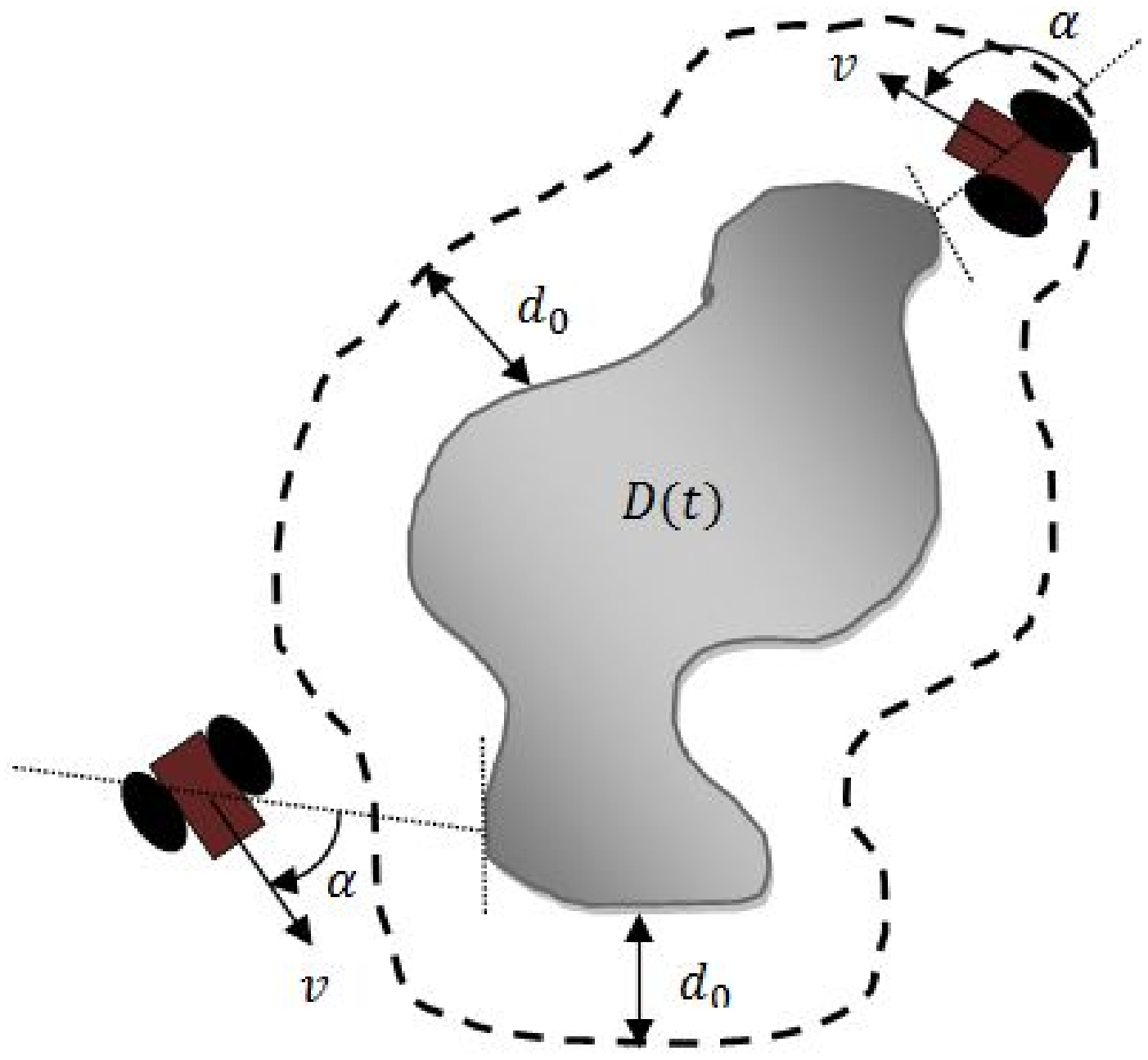}}
\label{sli}}
\hfill
\subfigure[]{\scalebox{0.25}{\includegraphics{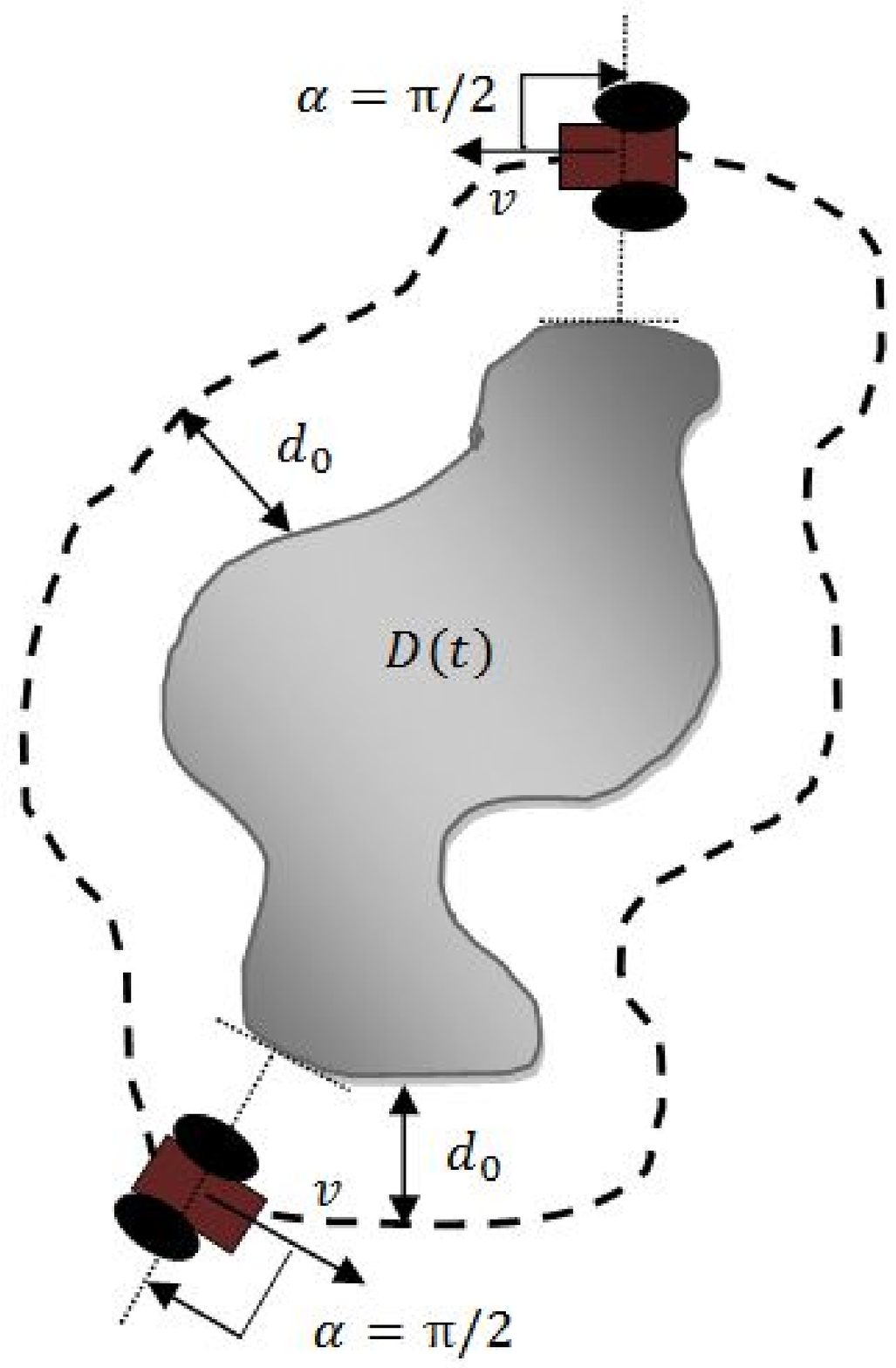}}
\label{lin}}
\caption{(a) Sliding motion towards the equidistant curve; (b) Motion along the equidistant curve}
\end{figure}
In more details, this motion looks as follows. The equation $\dot{d}+\chi(d-d_0)=0$ means that
the angle $\alpha$ between the relative velocity of the vehicle $\vec{v}_r$ and the line of sight at the nearest point of the obstacle equals $\alpha = \arccos\frac{\chi(d-d_0)}{v_r}$. Hence the angle $\alpha$ is obtuse for $d<d_0$ and acute for $d>d_0$, and so the vehicle is driven towards the desired distance $d_0$ to the obstacle; see Fig.~\ref{sli}. In doing so, $\alpha$ is kept constant $\alpha = \arccos[\frac{v_\ast}{v_r} \sgn \chi(d-d_0)]$ in the saturation zone $|d-d_0| > \delta$. As $d$ leaves this zone and approaches $d_0$, the angle goes to $\frac{\pi}{2}$.
In the limit where $d=d_0$, the vehicle is oriented parallel to the obstacle boundary $\alpha = \pi/2$, which means traveling along the $d_0$-equidistant curve; see Fig.~\ref{lin}.
\par
The control law \eqref{avoid} is activated in a close proximity of the obstacle. Otherwise the vehicle is driven towards the target in a straight line:
\begin{equation}
\label{pur}
u(t) \equiv 0, \qquad v(t) = \Upsilon[d(t)].
\end{equation}
Switching (\ref{pur}) $\mapsto$ (\ref{avoid}) occurs when the distance to the nearest obstacle does not exceed a given threshold $d_{\text{av}} \in (d_0, d_0^\Upsilon]$ and $\dot{d}+\chi(d-d_0) \leq 0$; switching (\ref{avoid}) $\mapsto$ (\ref{pur}) occurs when the vehicle is headed for the target and $\dot{d}+\chi(d-d_0) \geq 0$.
\par
The second relation from \eqref{pur} does not mean that the vehicle should have constant access to the distance $d(t)$: it suffices that this distance can be measured only in a vicinity of the obstacle $d(t) \leq C$, where $C\geq d_{\text{cr}}$.
\par
Practically the discrepancy $d_{\text{cr}}-d_{\text{av}}$ and the profile of $\Upsilon (\cdot)$ on $\big[d_{\text{av}},d_{\text{cr}}\big]$ are chosen with regard to the acceleration capacity of the wheelchair.
\section{Main assumptions}
\label{S3}
To describe the properties of the moving continuum $D(t)$,\footnote{Discussion of the properties of the moving continuum $D(t)$ basically follows that in \cite{MaChaSa12a}.} we introduce a {\it reference configuration} $D_\ast \subset \br^2$ of the body and the
{\it configuration map} $\Phi(\cdot,t): \br^2 \to \br^2$ that transforms $D_\ast$ into the current configuration $D(t) = \Phi[D_\ast,t]$.
We limit ourselves by only few conventions typical for the continuum mechanics and listed in the following.
\begin{assumption} \label{ass}
The reference domain $D_\ast$ is compact and has a smooth boundary $\partial D_\ast$. The configuration map $\Phi(\cdot,t)$ is defined on an open neighborhood of $D_\ast$ and is smooth, the determinant of its Jacobian matrix is everywhere positive $\det \Phi^\prime_{\blr} >0$, and $\Phi(\blr_1,t) \neq \Phi(\blr_2,t)$ whenever $\blr_1 \neq \blr_2$.
\end{assumption}
\par
Hence the boundary $\partial D(t)$ is smooth at any time.
\par
To proceed, we need some notations, which partly refer to the Eulerian formalism in description of the body motion:
\begin{itemize}
\item $\Phi^{-1}(\cdot,t)$ --- the inverse map;
\item $\vec{V}(\blr,t):= \Phi^\prime_t[\Phi^{-1}(\blr,t)]$ --- the velocity vector-field;
\item $\vec{A}(\blr,t):= \Phi^{\prime\prime}_{tt}[\Phi^{-1}(\blr,t)]$ --- the acceleration field;
\item $ \vec{V}^\prime_{\blr}[\blr,t]$ --- the spatial velocity gradient tensor;
\item $\mathscr{E}:= \frac{1}{2} \left[ \vec{V}^\prime_{\blr}[\blr,t]+  \vec{V}^\prime_{\blr}[\blr,t]\vphantom{V}^{\top} \right]$ --- the strain rate tensor;
\item $\omega=\omega(\blr,t)$ --- the vorticity, or angular velocity of the rigid-body-rotation i.e., $\left( \begin{smallmatrix} 0 & \omega \\ - \omega & 0 \end{smallmatrix} \right) = \frac{1}{2} \left[ \vec{V}^\prime_{\blr} -  \vec{V}^\prime_{\blr}\vphantom{V}^{\top} \right]$;
    \item $\sclp{\cdot}{\cdot}$ --- the standard inner product;
    \item $\vec{T}(\blr,t),\vec{N}(\blr,t)$ --- the Frenet frame of $\partial D(t)$ at $\blr \in \partial D(t)$; the domain $D(t)$ is to the left to the tangent vector $\vec{T}$, the normal vector $\vec{N}$ is directed inwards $D(t)$;
   \item $W_{T,t}(\blr,t):= \langle\vec{W}(\blr,t); \vec{T}(\blr,t)\rangle$ --- the tangential component of the vector-field $\vec{W}$ at $\blr \in \partial D(t)$ at time $t$;
       \item $W_{N,t}(\blr,t)$ --- the normal component;
       \item $\kappa(\blr,t)$ --- the signed curvature of $\partial D(t)$ at the point $\blr$;
\item $\sigma [\blr,t]:= \langle\mathscr{E}\vec{T}; \vec{N} \rangle + \omega = \langle \vec{V}^\prime_{\blr}\vec{T};\vec{N}\rangle$ --- the normal component of the velocity rate-of-change under an infinitesimally small shift along the boundary $\partial D(t)$;
 \item $\blr(t)$ --- the current position of the vehicle;
 \item $\blr_\ast(t)$ --- the point of $\partial D(t)$ nearest to $\blr(t)$.
\end{itemize}
The signed curvature $\kappa$ assumes positive and negative values on the convex and concave parts of $\partial D(t)$, respectively.
\par
The next assumption is typically fulfilled in real world, where physical quantities take bounded values.
\begin{assumption}
\label{bounded}
The scalars $V_N[\blr,t]$, $V_T[\blr,t]$, $A_N[\blr,t]$, $\sigma[\blr,t]$, $\kappa[\blr,t]$ remain bounded as $\blr \in \partial D(t)$ and $t \to \infty$.
\end{assumption}
\par
As $t$ runs within any finite time horizon, these quantities do remain bounded over $\blr \in \partial D(t)$ since all of them, along with $\Phi(\blr,t)$, continuously depend on $\blr$ and $t$, and the reference domain $D_\ast$ is compact by Assumption~\ref{ass}.
\par
In order that the vehicle be capable of bypassing the obstacle at the desired distance $d_0$, the $d_0$-equidistant curve to be traced should not be sharply contorted and move too fast.
The first and second parts of these requirements are fleshed out by
the following Assumption~\ref{ass.focus}
and Lemma~\ref{lemma.2}. The former is shown to be nearly unavoidable even in the case of the static $D(t) \equiv D$ domain \cite{MTS11}, whereas the proof of the latter can be found in \cite{MaChaSa12a}.
\begin{assumption}
\label{ass.focus} For any point $\blr$ on the $d_0$-equidistant curve of $D(t)$, the distance from $\blr$ to $D(t)$ is furnished by a single point $\blr_\ast \in D(t)$, and $1+ d_0 \kappa(\blr_\ast,t) > 0$. Moreover,
\begin{equation}
\label{positive}
\varliminf_{t \to \infty}\inf_{\blr_\ast \in \partial D(t)}
\left( 1+ d_0 \kappa[\blr_\ast,t] \right)>0.
\end{equation}
\end{assumption}
\begin{lemma}\label{lemma.2}
Let the vehicle travel so that it constantly overtakes the nearest boundary point of the obstacle and $d(t) \equiv d_0, v(t) \equiv v$.
 For any $t$, the parameters of the obstacle motion at the point $[\blr_\ast(t),t]$ satisfy the following relations:
\begin{multline}
\left| V_{N}\right| \leq v, \quad  \frac{|\mathscr{A}|L}{\sqrt{v^2-V_N^2}} + v \leq V , \\ \sqrt{v^2 - V_N^2} \geq \pm \left( V_T + d_0 \sigma \right),
\\
\mathscr{A}= \mathscr{A} (\blr_\ast,t,d_0,v):=A_N+ \frac{2 \sigma \xi + \kappa \xi^2 - d_0 \sigma^2}{1+\kappa d_0}\\
 \text{\rm and} \qquad \xi := -V_T \pm \sqrt{v^2 - V^2_N}.
\label{basic.nonstrict}
\end{multline}
Here the sign $+$ is taken if the vehicle moves so that the obstacle is to the left, and $-$ is taken otherwise.
\end{lemma}
\par
The assumption $v(t) \equiv v$ follows from $d(t) \equiv d_0$ if the speed is a function of the distance $d(t)$, like in \eqref{avoid} and \eqref{pur}.
Thus conditions \eqref{basic.nonstrict} must be satisfied, otherwise
the control objective cannot be achieved. Since the nearest boundary point $\blr_\ast(t)$ is not known in advance, it is reasonable to require that \eqref{basic.nonstrict} holds for all boundary points at any time.
The next assumption enhances a bit this by substituting the uniformly strict inequality in place of the non-strict one.
\begin{assumption}\label{ass.1}
There exist $\lambda_v, \lambda_a \in (0,1), \varepsilon_v>0$, $v_0 \in (0,V)$ such that for $v:=v_0$, the following inequalities hold:
\begin{eqnarray}
\label{ass.ineq}
\left| V_{N}\right| \leq \lambda_v v, \quad  \frac{|\mathscr{A}|L}{\sqrt{v^2-V_N^2}} + v \leq \lambda_a V,
\\
\sqrt{v^2 - V_N^2} \geq \left| V_T + d_0 \sigma \right| + \varepsilon_v
\label{ass.ineq1}
\end{eqnarray}
at any time $t$, point $\blr_\ast \in \partial D(t)$, and with both signs $\pm$.
\end{assumption}
\par
For any finite time horizon, such $\lambda_v,\lambda_a,\varepsilon_v$ exist if and only if the inequalities in \eqref{ass.ineq} are strict. Assumption~\ref{ass.1}  ensures that the strict inequalities do not degrade as $t \to \infty$.
This assumption clearly holds for motionless obstacles ($A_N=0, V_T=V_N =0, \sigma=0, \xi = \pm v$) with $v=v_0 \approx 0$.
\par
Inequalities \eqref{ass.ineq} allow not only to maintain the distance to the obstacle $\dot{d}=0$ but also to both decrease $\dot{d}<0$ and increase $\dot{d}>0$ it via proper manipulation of the control $u$, thus making the output $d$ locally controllable. As may be shown, the above assumptions are enough for local stability of DDWR in a vicinity of the desired equidistant trajectory. To establish non-local convergence to this trajectory, we need to ensure controllability during the transient.
To this end, we define the {\it launching motion} as that with the speed $v_0$ and maximal turning rate $u \equiv - (V-v_0)/L$ starting at a time $t_\ast$ from a position at the distance $d_{\text{av}}$ from the obstacle $D(t_\ast)$. This motion is over a circle of the radius $v_0L/(V-v_0)$ and is nothing but the motion that is typically observed at the initial stage of obstacle avoidance maneuver.
\par
Our final assumption encompasses the previous nearly unavoidable assumptions.
\begin{assumption}
\label{ass.main}
There exist $v_0 \in (0,V)$, $\lambda_v, \lambda_a \in (0,1), \varepsilon_v>0$, $d_-<d_+$ such that the following claims hold:
\begin{enumerate}
\item $d_{\text{safe}} < d_- < d_0 < d_+$ and \eqref{positive}---\eqref{ass.ineq1} are true for $v:=v_0, d_0:=d$ and any $d \in [d_-,d_+]$;
    \item During the first one and a half turn of any launching motion
    the distance $d(t)$ from DDWR to $D(t)$ lies in the interval $[d_-,d_+]$;
    \item In this motion, the total rotation angle $\alpha$ of $\blr(t) - \blr_\ast(t)$ does not exceed $(t_\ast-\tau_{\text{turn}})(V-v_0)/L$ at some time $t_\ast \in [\tau_{\text{turn}},\tau_{1.5 \text{turn}}]$, where $\tau_{\text{turn}}$ and $\tau_{1.5 \text{turn}}$ are the times when one full turn and {\rm 1.5} turns are completed, respectively.
\end{enumerate}
\end{assumption}
\par
\begin{figure}
\centering
\subfigure[]{\scalebox{0.3}{\includegraphics{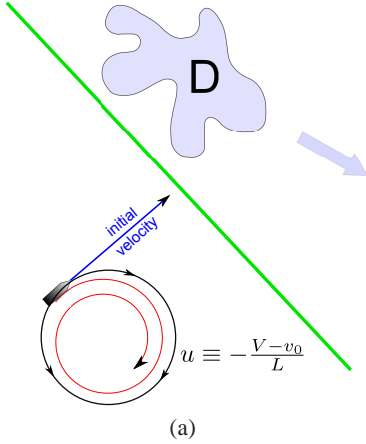}}\label{fig2a}}
\caption{Launching motion}
\end{figure}

 Here 3) holds with $t_\ast=\tau_{\text{turn}}$ for steady obstacles since then $\alpha =0$. For moving and deforming obstacles, 3) holds with $t_\ast= \tau_{\text{1.5 turn}}$ if e.g., the obstacle $D(t)$ and WC remain in disjoint steady half-planes and so $\alpha \leq \pi$; see Fig.~\ref{fig2a}.
\par
Since in \eqref{basic.nonstrict}, $\mathscr{A}$ is the (first degree)/(first degree) rational function in $d_0$, simple calculus show that it suffices to check the second inequality from \eqref{ass.ineq} only for $d = d_-,d_+$. Similarly \eqref{positive} holds for all $d \in [d_-,d_+]$ if and only if
$$
\varliminf_{t \to \infty}\inf_{\blr \in \partial D(t): \kappa[\blr,t] <0}
\left( 1+ d_+ \kappa[\blr,t] \right)>0.
$$
\section{Tuning controller parameters and the main results}
\label{sec.tmr}
In this section, we consider the parameters $\lambda_v,\lambda_a <1$, $d_-<d_+$, $v_0$ taken from Assumption~\ref{ass.main}. To tune the controller parameters, we first pick
two reals  $\eta_v >0, \eta_a >0 $ so that
\begin{equation}
\label{lambda-eta}
\lambda_v+\eta_v <1,\qquad \lambda_a+\eta_a <1.
\end{equation}
For any choice of the sign in $\pm$ and $v:=v_0$, the function
\begin{multline}
\label{func}
\Omega(\blr_\ast,t,d,z):=L \left| \frac{A_N+ \frac{2 \sigma \ov{\xi} + \kappa \ov{\xi}^2 - d \sigma^2}{1+\kappa d}}{\sqrt{v^2 - (V_N+z)^2}} \right| +v,\\
\text{where}\quad \ov{\xi} := -V_T \pm \sqrt{v^2 - (V_N+z)^2}
\end{multline}
is continuous in $z \approx 0$ uniformly over $\blr_\ast \in \partial D(t), t \geq 0, d \in [d_-,d_+]$ and for $z:=0$, equals the left-hand side of the second inequality from \eqref{ass.ineq}. Hence there exists $z_\ast >0$ such that for above $\blr_\ast,t,d$,
\begin{equation}
\label{eneq}
\Omega(\blr_\ast,t,d,z) < (\lambda_a + \eta_a) V \qquad \forall z \in [-z_\ast,z_\ast].
\end{equation}
\par
Now we are in a position to state the main theoretical result.
\begin{theorem}
\label{th.m} Let Assumptions {\rm \ref{ass}, \ref{bounded}}, and {\rm \ref{ass.main}} be satisfied, the parameter $v_0$ of the controller be taken from Assumption~{\rm \ref{ass.main}}, and the other parameters be chosen so that $v_{\text{cr}} \in (0,V)$,
\begin{multline}
\label{main.cond1}
v_\ast = \gamma \delta \leq \min\{\eta_v v_0, z_\ast \}, \\
(\lambda_a + \eta_a)  + \frac{\gamma L v_\ast}{v_0(V-v_0) \sqrt{1-(\lambda_v+\eta_v)^2}} <1 .
\end{multline}
Suppose that after the distance from DDWR to $D(t)$ reduces to the dangerous level $d_{\text{av}}$, the avoidance control law \eqref{avoid} is activated and then drives the vehicle.
Then the distance to the obstacle constantly exceeds the safety margin $d_{\text{\rm safe}}$ and moreover, lies within the interval $[d_-,d_+]$ from Assumption~{\rm \ref{ass.main}} and goes to the desired value $d(t)\xrightarrow{t \to \infty} d_0, \dot{d}(t)\xrightarrow{t \to \infty} 0$.
Since some time instant, DDWR constantly overtakes the nearest point of the obstacle boundary.
\end{theorem}
The proof of this theorem is given in Section~\ref{sec.proof}.
\par
Chosen $\lambda_{\nu}, \eta_\nu, \nu=v,a$ and $\gamma$,
the parameter requirements \eqref{main.cond1} can be always satisfied by picking $\delta$ small enough.
\begin{remark}
\rm
\label{rem.reverse}
Theorem~\ref{th.m} remains true for the navigation law \eqref{avoid} with the reversed sign
\begin{equation}
\label{c.a.rev}
\begin{array}{l}
	u(t) = - \frac{V - v(t)}{L} \cdot \sgn \big\{ \dot{d}(t) + \chi[d(t)-d_0]\big\}
\\
v(t) = \Upsilon[d(t)]
\end{array}
 \end{equation}
provided that in the definition of the launching motion, the turning rate $u \equiv (V-v_0)/L$ is put in place of $u \equiv - (V-v_0)/L$.
\end{remark}
\par
Target reaching property requires much more sophisticated analysis even in the case of steady obstacles.
In general, the control law \eqref{avoid} does not ensure that dealing with a given obstacle will terminate \cite{MaHoSa11ar}. However, this may hold only for complex maze-like obstacles \cite{MaHoSa11ar}.
Even in this case, there is a remedy in the form of a random choice among the control laws \eqref{avoid} and \eqref{c.a.rev} whenever a new avoidance maneuver is commenced \cite{MaHoSa11ar}. Unavoidable target reaching in the scene with steady convex obstacles was theoretically demonstrated in \cite{MTS11} for a unicycle driven by a similar control law.

\section{Proofs of theoretical results}
\label{sec.proof}
From now on, the Frenet frame $[\vec{T}(\blr,t), \vec{N}(\blr,t)]$ and the variables attributed to the motion of the domain, like $\vec{V}(\blr,t), \vec{A}(\blr,t)$, etc., are considered only for $\blr:=\blr_\ast(t)$. With a slight abuse of notations, the resultant argument $[\blr_\ast(t),t]$ is replaced by $t$ and dropped if  $t$ is clear from the context.
\begin{lemma}[\cite{MaChaSa12a}]
The velocity $\vec{v}$ of the robot has the form
\begin{equation}
\label{vel.0}
\vec{v}= \underbrace{\big[ \overset{\bul}{s}+d\mu\big]}_{\xi}\vec{T} + \vec{V} - \dot{d}\vec{N},
\end{equation}
where $\overset{\bul}{s}$ is the speed of the relative motion of $\blr_\ast$ along the boundary $\partial D$, i,e, $\dot{\blr}_\ast -\vec{V} = \overset{\bul}{s} \vec{T}$,
and
\begin{equation}
\label{mmuu}
\mu := \left\langle \vec{V}^\prime_{\blr} \vec{T}, \vec{N} \right\rangle + \kappa\overset{\bul}{s} =
\left\langle \mathscr{E} \vec{T}, \vec{N} \right\rangle +\omega+ \kappa\overset{\bul}{s}.
\end{equation}
\end{lemma}
\begin{lemma}[\cite{MaChaSa12a}]
The following relations hold
\begin{gather}
\ddot{d}
 =
\xi \mu - u \left[ \xi + V_T\right] + A_N
\nonumber
\\
+ \overset{\bul}{s} \sclp{\vec{V}^\prime_{\blr}\vec{T}}{\vec{N}}
\label{dddd}
 -\frac{\dot{v}}{v} (V_N - \dot{d})
\\
= - u \left[ \xi + V_T\right] + A_N+ \frac{2 \sigma \xi + \kappa \xi^2 - d \sigma^2}{1+\kappa d}
\nonumber
\\
\label{dddd1}
-\frac{\dot{v}}{v} (V_N - \dot{d});
\\
\overset{\bul}{s} = \frac{\xi - d \sigma}{1+\kappa d},
\label{s.dot}
\end{gather}
where $\overset{\bul}{s}$ is the speed of the relative motion of $\blr_\ast$ along the boundary $\partial D$, i,e, $\dot{\blr}_\ast -\vec{V} = \overset{\bul}{s} \vec{T}$.
\end{lemma}
 \par
 {\it Proof of Lemma}~\ref{lemma.2}. Since $d(t) \equiv d_0 \Rightarrow \dot{d}(t) \equiv 0, \ddot{d}(t) \equiv 0, \dot{v}(t) = F^\prime[d(t)] \dot{d}(t) \equiv 0$,
 \eqref{vel.0} yield that $v_N = V_N, v_T = \xi +V_T$, where $v_T^2+v_N^2 = v^2 \leq V^2$, which implies the first and last relations from \eqref{basic.nonstrict}.
 Similarly we have by \eqref{dddd1},
$$
 u  = \pm \frac{1}{\sqrt{v^2-V_N^2}}\left[ A_N+ \frac{2 \sigma \xi + \kappa \xi^2 - d \sigma^2}{1+\kappa d}\right].
$$
The proof is completed by noting that
$|u| \leq \frac{V-v}{L}$ by \eqref{control.constr}, and overtaking means that $\overset{\bullet}{s}>0$ if the vehicle has the domain to the
left, and that $\overset{\bullet}{s}<0$ otherwise.
\par
From now on, the assumptions of Theorem~\ref{th.m} are assumed to hold and we consider the wheelchair driven by \eqref{avoid}.
\begin{lemma}
\label{lem.cor1} Within the domain $d \in [d_- ,d_+]$, the surface
$S:=\dot{d} + \chi(d-d_0)=0$ is sliding in the sub-domain
\begin{equation}
\label{main.cond2} \xi+V_T >0
\end{equation}
and two-side repelling in the sub-domain
\begin{equation}
\label{d.asta} \xi+V_T <0
\end{equation}
if \eqref{main.cond1} holds. On this surface within the above
domain,
\begin{equation}
\label{in.sin} |\xi+V_T| \geq v \sqrt{1-(\lambda_v+\eta_v)^2}>0.
\end{equation}
\end{lemma}
\pf
Whenever $S:=\dot{d}-\nu=0$, where $\nu:= - \chi(d-d_0)$, we have by \eqref{chi},
\begin{equation}
\label{nu.ets} |\dot{d}|=|\nu| \leq v_\ast; \quad |\dot{\nu}| \leq \gamma |\dot{d}| = \gamma |\chi(d-d_0)| \leq \gamma
v_\ast .
\end{equation}
Due to \eqref{vel.0},
\begin{equation}
\label{v.sq}
(\xi+V_T)^2 + (V_N-\dot{d})^2 = v^2.
\end{equation}
and so
\begin{gather}
\nonumber
|\xi+V_T| = \sqrt{v^2 - (V_N-\dot{d})^2} \geq \sqrt{v^2 - (|V_N|+v_\ast)^2}
\\
\label{implic1}
\overset{\text{(a)}}{\geq} v \sqrt{1 - (\lambda_v+\eta_v)^2} \Rightarrow \text{\eqref{in.sin}},
\end{gather}
where (a) holds due to \eqref{ass.ineq} and \eqref{main.cond1}.
Furthermore,
\begin{gather}
\nonumber
\dot{S}=\frac{d}{dt} \big[ \dot{d} + \chi(d-d_0) \big] = \ddot{d} -
\dot{\nu}
\\
\nonumber
\overset{\text{\eqref{dddd1}}}{=\!=} - u \left[ \xi + V_T\right] + A_N+ \frac{2 \sigma \xi + \kappa \xi^2 - d \sigma^2}{1+\kappa d}  - \dot{\nu} ,
\end{gather}
where $\xi = -V_T \pm \sqrt{v^2-(V_N-\dot{d})^2}$ due to the starting argument from  \eqref{implic1}.
Let $B(\cdot)$ stand for the expression within $|\cdot|$ in \eqref{func} and $\ov{u}:= (V-v_0)/L$. Then we see that
\begin{gather}
\dot{S} = \left[ \xi + V_T\right] \left[ -u + B(\blr_\ast,t,d,\dot{d}) - \frac{\dot{\nu}}{\xi + V_T}\right]
\end{gather}
Here due to \eqref{eneq},
\begin{multline}
\label{omega.est}
\left|B(\blr_\ast,t,d,\dot{d}) - \frac{\dot{\nu}}{\xi + V_T}\right|
\leq (\lambda_a + \eta_a) \ov{u} + \frac{\gamma v_\ast}{|\xi + V_T|}
\\
\overset{\text{\eqref{in.sin}}}{\leq} (\lambda_a + \eta_a) \ov{u} + \frac{\gamma v_\ast}{v \sqrt{1-(\lambda_v+\eta_v)^2}} \overset{\text{\eqref{main.cond1}}}{<} \ov{u}
\end{multline}
So the signs
taken by $\dot{S}$ for $u=\ov{u}$ and $u=-\ov{u}$, respectively, are
opposite. In the case \eqref{main.cond2}, the sign is opposite to $\sgn u \overset{\text{\eqref{avoid}}}{=} \sgn S$; in the case \eqref{d.asta} the signs are equal. This implies the conclusion of
the lemma.
\epf
\begin{lemma}
\label{lem.attr}
If the equation $\dot{d} +
\chi(d-d_0)=0$ becomes true at some time $t_0$ when $d \in
[d_-,d_+]$, then monotonically and exponentially fast $d \to d_0, \dot{d} \to 0$ as $t \to \infty$. Furthermore, $d(t) \in [d_- ,d_+]$ and \eqref{main.cond2} holds for all $t \geq t_0$.
\end{lemma}
\pf Lemma~\ref{lem.cor1} guarantees that first, $\xi+V_T>0$ at
$t=t_0$ and second, this inequality is still valid and sliding
motion occurs while $d \in [d_-,d_+]$. During this motion, $\dot{y}
=- \chi(y)$ for $y:= d-d_0$, where $\chi(y)\cdot y>0 \; \forall y
\neq 0$ and $\chi(0)=0$. It follows that any solution $d$ of the
sliding mode differential equation monotonically converges to $d_0$.
At the same time, $d_0 \in [d_-,d_+]$.
Hence $d$ will never leave the interval $[d_-,d_+]$, the sliding
mode will never be terminated, the inequality $\xi+V_T>0$ will never be violated, and $d \to d_0, \dot{d} \to 0$ as $t \to \infty$.
Application of the Lyapunov's first method to the equation $\dot{y}
=- \chi(y)$ at the equilibrium point $y=0$ shows that the
convergence is exponentially fast.
\epf
\begin{lemma}
\label{lem.fin} Both relations $\dot{d} + \chi(d-d_0)=0, \xi+V_T>0$ become true at some moment following that when the distance $d$ reduces to $d_{\text{av}}$. At this moment, $d \in [d_-,d_+]$ and moreover, $d \geq d_-$ until this moment.
\end{lemma}
\pf
Let $t=t_0$ be the time when the distance reduces to $d_{\text{av}}$.
If the required relations hold at $t=t_0$, the claim is true.
Otherwise, $\dot{d} + \chi(d-d_0)\neq 0$ for $t>t_0, t \approx t_0$.
\par
Suppose that $\dot{d} + \chi(d-d_0) > 0$ for $t>t_0, t \approx t_0$. Since the control law \eqref{avoid} is ultimately activated by the assumptions of Theorem~\ref{th.m},
$\dot{d} + \chi(d-d_0)$ arrives at zero at some time $t_1>t_0$ such that $d(t_1) \leq d_{\text{av}}$. By the argument following \eqref{pur}, $d(t) \geq d_0 \geq d_- \;\forall t \in [t_0,t_1]$, and $d(t_1) \in [d_0,d_{\text{av}}] \subset [d_-,d_+]$. It remains to note that $\xi+V_T<0$ is impossible by Lemma~\ref{lem.cor1}.
\par
Suppose that $\dot{d} + \chi(d-d_0) < 0$ for $t>t_0, t \approx t_0$.
Since $t=t_0$ and
until the first time $t_1$ when the equation $\dot{d} + \chi(d-d_0) = 0$ becomes true, the vehicle moves with the constant control $u \equiv - \ov{u}, \ov{u} := (V-v_0)/L$.
Now we analyze the motion of the robot driven by the constant control $u \equiv
-\ov{u}$ for $t \in [t_0,t_0+\tau_{\text{1.5 turn}}]$. By definition, this is a launching motion  and the vehicle velocity $\vec{v}(t)$ rotates clockwise at the rate $\ov{u}$. By the last claim from Assumption~\ref{ass.main}, there exists $t_\ast \in \big[t_0+\tau_{\text{turn}}, t_0+\tau_{\text{1.5 turn}}\big]$ such that as $t$ runs from $t_0$ to $t_\ast$,
the total rotation angle of the vector $\blr(t) -\blr_\ast(t)$ does not exceed
$\beta:=(t_\ast-\tau_{\text{turn}} -t_0)\ov{u}$. Meanwhile $\vec{v}(t)$ rotates through the angle $2 \pi + \beta$.
Hence there exist two time instants $s_i \in \left[t_0, t_0+\tau_{\text{1.5 turn}}\right], i=1,2$ such that $\vec{v}(s_i)$ and $(-1)^i [\blr(s_i) -\blr_\ast(s_i)]$ are co-linear and identically directed for $i=1,2$. Thanks to \eqref{vel.0},
\begin{gather*}
(-1)^i\dot{d}(s_i) = v+ (-1)^i V_N(s_i) \geq v - |V_N| \overset{\text{\eqref{ass.ineq}}}{\geq} (1-\lambda_v) v
\\
\overset{\text{\eqref{lambda-eta}}}{>} \eta_v v \overset{\text{\eqref{main.cond1}}}{\geq} v_\ast \geq |\chi[d(s_i) - d_0]|
\end{gather*}
Thus the continuous function of time $S$ assumes values of opposite signs at $t=s_1,s_2$. It follows that $S$ inevitably arrives at zero within $\left[t_0, t_0+\tau_{\text{1.5 turn}}\right]$. Since $d \in [d_- ,d_+]$ at this moment by (2) in Assumption~\ref{ass.main}, Lemma~\ref{lem.cor1} implies that
\eqref{d.asta} cannot be true at this moment. So \eqref{main.cond2} is true due to \eqref{in.sin}, which and completes the proof.
\epf
\par
{\it Proof of Theorem}~\ref{th.m}. This theorem is immediate from
Lemmas~\ref{lem.attr}---\ref{lem.fin}, along with \eqref{s.dot} and \eqref{implic1}.
The last claim of the theorem is straightforward from \eqref{s.dot} and \eqref{v.sq}.

\bibliographystyle{IEEEtr}
\bibliography{wheelchair,Hamidref}
\end{document}